\RequirePackage[l2tabu, orthodox]{nag}

\documentclass[a4paper,reqno]{amsart}

\usepackage{lmodern}
\usepackage[T1]{fontenc}
\usepackage[utf8]{inputenc}
\usepackage[english]{babel}
\usepackage{microtype} 

\usepackage{amsmath,amssymb,amsthm,mathrsfs,latexsym,mathtools,mathdots,enumerate,tikz,ytableau}
\usetikzlibrary{positioning}
\usepackage[enableskew,vcentermath]{youngtab}
\usepackage{graphicx}
\usepackage{hyperref}
\usepackage{bm}
\usepackage{todonotes}

\usepackage{ifpdf}
\graphicspath{{./figures/}}
\usepackage{epstopdf}
\DeclareGraphicsExtensions{.png,.jpg,.eps,.epsf}

\def\la{\lambda}
\def\ka{\kappa}

\def\si{\sigma}

\renewcommand{\phi}{\varphi}

\title
{On the combinatorial local log-concavity conjecture and a result of Stanley}
\author{Valentin Féray}
\address{Institut für Mathematik, Universität Zürich, Winterthurerstrasse 190, 8057 Zürich, Switzerland}
\email{valentin.feray@math.uzh.ch}

\begin{document}

\maketitle

The purpose of this note is to explain that the {\em combinatorial local log-concavity conjecture} introduced by Gross, Mansour, Tucker and Wang
\cite[Conjecture 3.1]{CLLC}
in fact follows from a result of Stanley \cite[Corollary 3.3]{Stanley}.
\bigskip

We use the following notation:
\begin{itemize}
    \item $\la$ is a partition of $n$ and $\pi$ a fixed permutation of type $\la$.
    \item $Q_n$ is the set of $n$-cycles in the symmetric group $S_n$.
    \item $\rho(\si)$ denotes the cycle-type of $\rho$,
        while $\kappa(\si)$ denotes the number of cycles of a permutation $\si$.
    \item $z_\la$ is the usual combinatorial factor,
        that is the size of the centralizer of $\pi$ in $S_n$.
\end{itemize}
In his article \cite{Stanley}, Stanley considers the following polynomial
\[P_\la(q) =  \sum_{\rho(w)=\la} q^{\ka\big( (1,\dots,n) w \big)},\]
while, in \cite{CLLC}, the authors introduce
\[F_\la(q) = \sum_{\zeta \in Q_n} q^{\lfloor (\ka(\zeta \pi)-1)/2 \rfloor}.\]
The combinatorial local log-concavity conjecture
states that for any partition $\la$, the polynomial
$F_\la$ is log-concave.
\medskip

The polynomials $P_\la$ and $F_\la$ are in fact closely related,
as shown by the following computation
\begin{align*}
    P_\la(q) &= \sum_{\rho(w)=\la} q^{\ka\big( (1,\dots,n) w \big)}\\
    &= \frac{1}{z_\la} \sum_{\si \in S_n} q^{\ka\big( (1,\dots,n) \si \pi \si^{-1} \big)}\\
    &= \frac{1}{z_\la} \sum_{\si \in S_n} q^{\ka\big(\si^{-1}  (1,\dots,n) \si \pi \big)}\\
    &= \frac{n}{z_\la} \sum_{\zeta \in Q_n} q^{\ka\big(\zeta \pi \big)}.
\end{align*}
Now we have two cases: 
\begin{itemize}
    \item either $n+\kappa(\pi)$ is even, in which case $\ka\big(\zeta \pi \big)$
        is odd for any $n$-cycle $\zeta$ in $Q_n$ and 
        \[P_\la(q)=\frac{n}{z_\la} q F_\la(q^2) ;\]
    \item or $n+\kappa(\pi)$ is odd, in which case $\ka\big(\zeta \pi \big)$
        is even for any $n$-cycle $\zeta$ in $Q_n$ and 
        \[P_\la(q)=\frac{n}{z_\la} q^2 F_\la(q^2). \]
\end{itemize}
From \cite[Corollary 3.3]{Stanley} (see also the discussion below this corollary),
we know that $P_\la(q)$ has only purely imaginary root,
which implies in both cases that $F_\la$ has (nonnegative) real roots.
The log-concavity of its coefficients follows immediately.

\end{document}